\documentclass[12pt]{amsart}
\usepackage{amsmath,amsthm,amsfonts,amssymb}
\begin{document} 
\newcommand{\tensor}{\otimes}
\newcommand{\A}{{\mathbb A}}
\newcommand{\B}{{\mathbb B}}
\newcommand{\C}{{\mathbb C}}
\newcommand{\F}{{\mathbb F}}
\newcommand{\Gs}{{}^{\sigma}}
\newcommand{\N}{{\mathbb N}}
\newcommand{\Q}{{\mathbb Q}}
\newcommand{\Z}{{\mathbb Z}}
\renewcommand{\P}{{\mathbb P}}
\newcommand{\R}{{\mathbb R}}
\newcommand{\rc}{\subset}
\newcommand{\Aut}{\mathop{Aut}}
\newcommand{\Gal}{\mathop{Gal}}
\newcommand{\rank}{\mathop{rank}}
\newcommand{\trace}{\mathop{tr}}
\newcommand{\dimc}{\mathop{dim}_{\C}}
\newcommand{\Lie}{\mathop{Lie}}
\newcommand{\Hom}{\mathop{Hom}}
\newcommand{\Auto}{\mathop{{\rm Aut}_{\mathcal O}}}
\renewcommand{\Vec}{\mathop{{\rm Vec}}}
\newcommand{\alg}[1]{{\mathbf #1}}
\newtheorem*{definition}{Definition}
\newtheorem*{claim}{Claim}
\newtheorem{corollary}{Corollary}
\newtheorem*{Conjecture}{Conjecture}
\newtheorem*{SpecAss}{Special Assumptions}
\newtheorem{example}{Example}
\newtheorem*{remark}{Remark}
\newtheorem*{observation}{Observation}
\newtheorem*{fact}{Fact}
\newtheorem*{remarks}{Remarks}
\newtheorem{lemma}{Lemma}
\newtheorem{proposition}{Proposition}
\newtheorem{theorem}{Theorem}
\title{%
Non-linearizable Actions of Commutative Reductive Groups.
}
\author {J\"org Winkelmann}
\begin{abstract}
We generalize a construction of Freudenburg and Moser-Jauslin in
order to obtain an example of a non-linearizable action of a
commutative reductive group on the affine space for every
field $k$ of characteristic zero which admits a quadratic
extension.
\end{abstract}
%
\address{%
J\"org Winkelmann \\
Mathematisches Institut\\
Universit\"at Bayreuth\\
Universit\"atsstra\ss e 30\\
D-95447 Bayreuth\\
Germany\\
}
\email{jwinkel@member.ams.org\newline\indent{\itshape Webpage: }%
http://btm8x5.mat.uni-bayreuth.de/\~{ }winkelmann/
}
\maketitle
\section{Introduction and Summary}
Let $G$ be an algebraic group acting on the affine space $\A^n$.
If $G$ is reductive it is not so easy to find an example of
such an action which is not obviously linearizable.
At a time it was even conjectured that at least over an algebraically
closed field $k$ every action of a reductive group is
linearizable (\cite{Kam}).

But then Schwarz found a non-linearizable such action
(\cite{S}). His
example is an action of the complex group $O_2(\C)$
which is a semi-direct product of $\Z/2\Z$ and the multiplicative
group $G_m(\C)\simeq\C^*$.
Related examples have been deduced, but results about
commutative reductive groups remain rare.
See \cite{K} for more about these topics.

Recently, Freudenberg and Moser-Jauslin constructed an example
of a non-linearizable real algebraic action of $S^1$ (\cite{FMJ}).
$S^1$ is a non-standard real form of the multiplicative group $G_m$.
In this note we will go through the arguments of Freudenberg
and Moser-Jauslin in order to verify that in their construction
the field of real numbers can be replaced by an arbitrary
field of characteristic zero with a non-trivial form of the
multiplicative group.
\begin{theorem}
Let $k$ be a field of characteristic zero, $\bar k$ an algebraic
closure and $H^0$ a $k$-group. We assume that $H^0$ is isomorphic
to the multiplicative group $G_m$ over $\bar k$, but not over
$k$. (Such $k$-groups are often called {\em non-split $k$-forms}
of the multiplicative group $G_m$.)

Then there is an action of $H^0$ on the affine space $\A^4$,
defined over $k$, which is not linearizable over $k$.
\end{theorem}

As usual, an action of a $k$-group $H^0$ on the affine space
$\A^n$ given by a morphism $\mu:H^0\times\A^n\to \A^n$
is called {\em linearizable over $k$} if and only if there
exists an automorphism $\phi$ of  $\A^n$ as a $k$-variety
such that $v\mapsto \phi(\mu(g,\phi^{-1}(v)))$ is a linear map
for every $g\in H^0$.

\begin{corollary}
Let $k$ be a field of characteristic zero which contains a
non-square $\alpha\in k^*\setminus\left(k^*\right)^2$.

Then there is a one-dimensional commutative connected reductive
$k$-group $G$ with a non-linearizable action on $\A^4$.
\end{corollary}
\begin{proof}
As discussed in \S2 below the existence of a quadratic field extension
implies the existence of a non-split $k$-form of the multiplicative
group $G_m$.
\end{proof}
\begin{corollary}
Let $k$ be one of the following:
\begin{enumerate}
\item the field $\R$ of real numbers,
\item the field $\Q_p$ of $p$-adic numbers,
\item a number field,
\item a finitely generated extension of $\Q$,
\item the field of rational functions of a positive-dimensional
variety defined over some field $k_0$ of characteristic zero,
\item the field of meromorphic functions of a complex manifold
which admits a non-constant meromorphic function, 
\item $\Q^{ab}$, i.e.~the maximal abelian field extension of $\Q$.
\end{enumerate}

Then there is a one-dimensional commutative connected reductive
$k$-group $G$ with a non-linearizable action on $\A^4$.
\end{corollary}
\begin{proof}
In view of the preceding results it suffices to show that each of these
fields admits  a non-trivial quadratic extension.

For $\R$ one takes $\C=\R[i]$. For $k=\Q^{ab}$ we recall that there
exists a number field $E$ with the quaternion group $\Gamma$ as Galois group
(e.g.~$\Q\left[\sqrt{(\sqrt 2+2)(\sqrt 3+3)}\right]$, see \cite{Mil}).
The quaternion group $\Gamma$ is two-step nilpotent. Hence its
commutator group is a non-trivial commutative subgroup $\Gamma'$.
Let $L$ denote the intermediate field corresponding to $\Gamma'$.
Then $E/L$ is a quadratic extension. Furthermore,
$L\subset \Q^{ab}$ (because $\Gal(L/\Q)\simeq \Gamma'$ is abelian)
and $E\not\subset\Q^{ab}$
(because $\Gal(E/\Q)\simeq\Gamma$ is not abelian). 
Thus we obtain a
quadratic extension $E\Q^{ab}/\Q^{ab}$.

For all the other fields $k$ listed above there is a discrete
valuation. A discrete valuation gives in particular a group
homomorphism $v$ from the multiplicative group $k^*$ to the additive
group $(\Z,+)$. Evidently an element $x\in k^*$ can not be
a square unless $v(x)$ is even. Therefore the existence of a discrete
valuation implies that there exists a non-trivial quadratic extension.
\end{proof}
Thus most well-known fields are either algebraically closed
or admit a quadratic extension.
But of course there are also fields which are neither algebraically
closed nor admit a quadratic extension. For instance, let
$\Gamma$ be the absolute Galois group of $\Q$ and
consider all group homomorphisms $\rho$ from $\Gamma$ to
$2$-groups (=finite groups whose order is a power of $2$).
Then the intersection of the kernels $\cap_\rho\ker\rho$ defines
a closed subgroup $H\subset\Gamma$ such that the corresponding
(infinite) extension field $L$ of $\Q$ is neither algebraically
closed nor admits any quadratic extension.

\section{$k$-forms}
Here we recall some basic facts about {$k$-forms}, see \cite{Se}, Ch.~3, \S1.

Let $k$ be a field of characteristic zero,
$K/k$ be a Galois extension and $G$ a $k$-group.
One is interested in classifying all $k$-groups $H$ which are
$K$-isomorphic to $G$. It is a standard result that these ``$k$-forms
of $G$''
are classified by the Galois cohomology set
\[
H^1(\Gal(K/k),\Aut{}_k(G)).
\]
Here we are interested in the special case where $G=G_m$ is the multiplicative
group. Then $\Aut_k(G)\simeq\Z/2\Z$, because $G_m$ admits only two
automorphisms: the identity and $z\mapsto\frac{1}{z}$.
Thus the action of $\Gal(K/k)$ on $\Aut_k(G_m)$ is necessarily trivial
and consequently
\[
H^1(\Gal(K/k),\Aut{}_k(G_m)) \ \simeq \ \Hom(\Gal(K/k),\Z/2\Z).
\]
Hence $k$-forms of $G_m$ are in one-to-one correspondence with
group homomorphisms from $\Gal(K/k)$ to $\Z/2\Z$ and therefore
in one-to-one correspondence to quadratic field extensions of $k$
which are contained in $K$.

In concrete terms these $k$-forms of $G_m$ 
can be described as follows:

Let $K/k$ be a quadratic extension of fields in characteristic zero.
Then there is a non-square $\alpha\in k$ such that $K=k[X]/(X^2-\alpha)$.
Each element $z\in K$ can be written in unique way as $z=x+ty$ where
$t$ is a fixed element with $t^2=\alpha$ and $x,y\in k$.
A quadratic extension is necessarily Galois and the non-trivial
element of $\Gal(K/k)$ acts by $z=x+ty\mapsto x-ty$.
Therefore the {\em Norm homomorphism} $N_{K/k}:K^*\to k^*$ is
given by $N_{K/k}(x+ty)=(x+ty)(x-ty)=x^2-\alpha y^2$.
The group of all elements $z\in K^*$ with $N_{K/k}(z)=1$
can now be identified with the $k$-rational points of the
following $k$-group:
\[
H=\left\{ 
\begin{pmatrix} x & \alpha y \\ y & x \\
\end{pmatrix}
: x^2-\alpha y^2 = 1
\right\}.
\]
In the special case $\C/\R$ this yields the real one-dimensional
compact torus $S^1$.

Let us explain the connection with Galois cohomology. 
Let $K/k$ be a quadratic extension and let $\sigma$ be the non-trivial
element of $\Gal(K/k)$. Then the $k$-form described above corresponds
to the cocycle mapping $\sigma$ to the automorphism $z\mapsto z^{-1}$
of $G_m$.

More generally, if $X$ is a variety defined over $K$ with an
involution
$\tau$ defined over $k$, there is a $k$-form $Y$ of $X$ associated
to the cocycle $\sigma\mapsto\tau$ and there is a bijection
$Y(k)\simeq\{p\in X(K):{}^\sigma p=\tau(p)\}$.

\section{The criterion of Masuda-Petrie}
Let $G$ be a $k$-group and let $V,W$ be $G$-modules and
let $\Vec_G(V,W)$ be the set of $G$-vector bundles over 
the $G$-variety $V$ such that
the fiber over $0_V$ is isomorphic to $W$.

Let us assume that there exists a subgroup $I\subset G$
such that $(V\oplus W)^I=V\oplus\{0\}$ where $(V\oplus W)^I$ denotes
the set of fixed points for the $I$-action on $W\oplus V$.

Let $E\in \Vec_G(V,W)$. Apriori, triviality of $E$ as $G$-vector bundle
is a stronger condition than linearizability of the $G$-action on
the total space: In the second condition all isomorphisms from $E$
to $V\oplus W$
as variety may be used to linearize the $G$-action on $E$ while
in the first condition only those isomorphisms are admitted which
preserve the vector bundle structure, i.e., are compatible with
the projection on $V$ and linear on each fiber.
However,
Masuda and Petrie showed in \cite{MP}
that under the above mentioned additional
assumption both conditions are equivalent.

Let us recall their arguments, for the convenience of the reader
as well as in order to convince us that this works over any base
field in characteristic zero.
Assume that there is a $G$-isomorphism of varieties $\phi:E\to V\oplus W$.
Then $E^I$ is mapped onto $V\oplus\{0\}$. It follows easily that
the zero-section of $E\to V$ equals 
 the fixed point set $E^I$.
Moreover it follows that $\phi$ induces an isomorphism between
the respective normal bundles. But the normal bundle for the zero-section
in a vector bundle is isomorphic to the given vector  bundle itself.
In this way we obtain a trivialization of $E$ as a $G$-vector bundle.

\section{The Schwarz action}
Let $G$ be the non-trivial semi-direct product of $\Z/2\Z=\{e,\tau\}$ with
the multiplicative group $G_m$.
Then $G$ is an algebraic group defined over $\Q$ (and therefore
defined
over every field $K$ of characteristic zero).

Schwarz introduced in \cite{S}
an action of $G$ on $\A_4$ which is defined over $\Q$
and is not linearizable over $\C$ (and therefore a fortiori not
linearizable over any field $k$ of characteristic zero).
In concrete terms it is given as follows:

We write $\Z/2\Z=\{e,\tau\}$
and define the action on $E=\A^4=\A^2\times\A^2$ by
\[
\lambda:\left( \begin{pmatrix}a \\ b \\
  \end{pmatrix},
  \begin{pmatrix}x \\ y \\
  \end{pmatrix}\right)
\ \longrightarrow\ 
\left( \begin{pmatrix}
\lambda^2a \\ \lambda^{-2}b \\
  \end{pmatrix},
  \begin{pmatrix}\lambda^{3}x \\ \lambda^{-3}y \\
  \end{pmatrix}\right)
\]
for $\lambda\in G_m(\C)=\C^*$
and
\[
\left( \begin{pmatrix}a \\ b \\
  \end{pmatrix},
  \begin{pmatrix}x \\ y \\
  \end{pmatrix}\right)
\in \A^2\times \A^2\simeq\A^4.
\]
We let $\tau $ act by
\[
\tau:\left( \begin{pmatrix}a \\ b \\
  \end{pmatrix},
  \begin{pmatrix}x \\ y \\
  \end{pmatrix}\right)
\ \longrightarrow\ 
\left( \begin{pmatrix}b \\ a \\
  \end{pmatrix},
\begin{pmatrix}
1 + ab + (ab)^2 & -b^3 &  \\
 a ^3 & 1-ab  \\
\end{pmatrix}
\cdot\begin{pmatrix}y \\ x \\
  \end{pmatrix}\right).
\]

Moser-Jauslin showed that the action of $\Z/2\Z=\{e,\tau\}$
can be linearized over $\Q$ (see \cite{MJ}). In concrete terms
this linearization as obtained in \cite{MJ} can be described
as follows:
$\phi\circ\tau\circ \phi^{-1}$ is a diagonal linear endomorphism
where $\phi$
is the automorphism given as:
\[
\phi:
\left(
\begin{pmatrix}
a \\ b\\
\end{pmatrix}
,
\begin{pmatrix}
x \\ y\\
\end{pmatrix}
\right)
\longrightarrow
\left(
\begin{pmatrix}
a \\ b\\
\end{pmatrix}
,
\begin{pmatrix}
c_{11} & c_{12} \\
c_{21} & c_{22} \\
\end{pmatrix}
\cdot
\begin{pmatrix}
x\\
y\\
\end{pmatrix}
\right)
\]
with
\begin{align*}
c_{11}&= 2(1+a-2b-a^2b+2ab^2-b^3)\\
c_{12}&=2(1-2a+b+ab+a^4-2a^3b+a^2b^2)\\
c_{21}&= -2 -a -b +ab -b^2 +a^2b -b^3 \\
c_{22}&= 2+b+a +a^2 +ab-a^3 +a^2b-a^4+a^2b^2.
\end{align*}
Note that $\phi$ is defined over $\Q$.

Note furthermore that the map
\[
\left(
\begin{pmatrix}
a \\ b\\
\end{pmatrix}
,
\begin{pmatrix}
x \\ y\\
\end{pmatrix}
\right)
\longrightarrow
\begin{pmatrix}
a \\ b\\
\end{pmatrix}
\]
realizes $E$ as $G$-vector bundle of rank $2$ over $\A^2$.

\section{A non split form of the Schwarz action}

The algebraic group $G$ and the space $E=\A^4$ both admit an involution:
On $G$ we take conjugation by $\tau$ and on $E$ we take $\tau$ acting as
described in the preceding section.
If $\mu:G\times E\to E$ denotes the morphism defining the Schwarz action,
then 
\[
\mu(g^\tau,\tau v)=\tau\cdot g\cdot \tau^{-1}\cdot \tau v=\tau\mu(g,v).
\]
Thus this pair of involutions on $G$ and $E$ is compatible with the
Schwarz action. Furthermore, both involutions are defined over $\Q$.

We fix  a quadratic extension $K/k$ of fields of characteristic zero
and let $\sigma$ denote
the non-trivial element of $\Gal(K/k)$.
We fix an element $j\in K$ such that ${}^\sigma j=-j$.
Then $K=k\oplus jk$ and $K=k[j]\simeq k[X]/(X^2-\alpha)$
where $\alpha=j^2$.
Now we can obtain $k$-forms of $G$, $E$ and the action $\mu$
associated to the Galois cocycle which maps the non-trivial element
$\sigma$ of $\Gal(K/k)$ to the respective involution chosen above.
We call these $k$-forms $H$, $E_0$ and $\mu_0$ respectively.

In concrete terms this works as follows:

\[
H(k)=\left\{g\in G(K):{}^\sigma g=\tau g \tau^{-1}\right\}.
\]

Since the involution $\tau$ is defined over $\Q$, it commutes
with the Galois action of any Galois group for any Galois extension
$K/k$. 
We define
\[
E_0(k)=\{v\in E(K)\simeq K^4: 
\Gs v=\tau(v)\}
\]
Using the linearization $\phi$ of the involution $\tau$,
we see that
a vector $(z_1,\ldots,z_4)\in K^4\simeq\phi(E)(K)$ 
is contained in $\phi(E_0(k))$ iff
$jz_1,z_2,z_3,z_4\in k$.
Therefore $E_0(k)\simeq k^4$ and $E_0\tensor_kK=E$.
Now $H(k)$ is generated by $\tau$ and those $\lambda\in G_m(K)=K^*$
for which $\Gs\lambda\lambda=1$.
Since $\tau$ commutes with $\sigma\circ\tau$ and $E_0(k)$ is the
fixed point set of $\sigma\circ\tau$, the set $E_0(k)$ is stabilized
by $\tau$. Next consider $\lambda\in (H(k))^0$ and $v\in E_0(k)$. Then 
$\Gs\lambda\lambda=1$ and 
\[
\sigma\circ\tau(\lambda v)=
\tau\circ\sigma(\lambda v)=
\tau(\Gs\lambda \Gs v)=
\tau(\lambda^{-1}\tau v)=
\lambda v.
\]
Hence $\lambda v\in E_0(k)$. Thus $H(k)$ stabilizes $E_0(k)$.

We underline that for the Schwarz action of $G$ on $\A^4$ there
is one set of coordinates for which the connected component
$G^0$ containing the neutral element acts linearly and a second,
different set of coordinates for which $\tau$ acts linearly,
but no coordinate set for which the whole group $G$ acts linearly.

The $k$-form $E_0$ of $E=\A^4$ is isomorphic to $\A^4$ as a $k$-variety.
We deduced this using the set of coordinates in which $\tau$ is
linear. In these coordinates the $G^0$-action is not linear.
Therefore the $H^0$-action on $E_0$ need not be linearizable
(over $k$)
although the $G^0$-action on $E$ is.

We also remark that there is a unique fixed point $0$ 
for $H^0$ on $E_0$.
It follows that the $H^0$-action in $E_0$ is either not linearizable
at all or is isomorphic to the induced $H^0$-action on the tangent space
$T_0E_0$ of this fixed point.

Furthermore, the involution $\tau$ as well as the automorphism of $E$ 
linearizing $\tau$ are compatible with the structure of $E$ as a 
$G$-vector bundle, which implies that $E_0$ inherits this structure,
i.e. is an $H$-vector bundle of rank $2$ over $\A^2$.

\section{Reduction from $H$ to $H^0$}
\newcommand{\cj}{\,{}^{\sigma}\!}
\begin{proposition}\label{prop1}
Let $K/k$ be a quadratic extension of fields of characteristic zero,
$\Gal(K/k)=\{id_K,\sigma\}$,
 $S=\{x\in
K^*:N_{K/k}(x)=x\cj x=1\}$
and $j\in K\setminus k$. 
For $\lambda\in S$ define
a morphism $\mu_\lambda:\A^2\to\A^2$ by 
\[
\mu_\lambda(x,y)=\left(\lambda^2x,\lambda^3y\right)
\]
and let $\tau:K^2\to K^2$ be a map with the following
properties:
\begin{enumerate}
\item[(i)]
$\tau(x,y)=\left(\cj x,\phi(x,y)\right)\ \forall x\in K$
for some map $\phi:K^2\to K$.
\item[(ii)]
The map $y\mapsto\phi(x,y)$ is $k$-linear in $y$ for every fixed 
$x\in K$.
\item[(iii)]
There is a $k$-morphism $\tau_0:\A^4\to\A^4$ such that
\[ 
\tau_0(z,w,u,v)=(a,b,c,d)
\]
whenever
\[
\tau(z+jw,u+jv) = (a+jb, c+jd)\quad (z,w,u,v,a,b,c,d\in k)
\]
\item[(iv)]
$\tau\circ\tau=id_{K^2}$.
\item[(v)]
$\tau\circ\mu_\lambda=\mu_{\frac{1}{\lambda}}\circ\tau$
for all $\lambda\in S$.
\end{enumerate}

\smallskip

Then there exists an element $\alpha\in S$
such that $\tau(x,y)=(\cj x,\alpha\cj y)$ for all $x,y\in K$.

In particular, the map $\tau:K^2\to K^2$ is $k$-linear.
\end{proposition}
\begin{remark}
In the formulation of the above proposition $S$ is a group
acting by $k$-linear transformations on $K^2\simeq k^4$. The group
$S$ and the map $\tau$ give a semi-direct product of $\Z/2\Z$ and $S$
acting on $K^2\simeq k^4$ 
and the proposition says that under certain assumptions the action of
this semi-direct product on $K^2\simeq k^4$ is necessarily also
by $k$-linear transformations.
\end{remark}

\begin{proof}
Condition $(iv)$ implies that $\tau$ is bijective.

Any $k$-linear map between $K$ vector spaces can be written as a
sum
of a $K$-linear and a $K$-antilinear map. Therefore condition $(ii)$
implies that there are maps $\phi',\phi'':K\to K$ such that
\[
\phi(x,y)=\phi'(x)y+\phi''(x)\cj y
\]
 for all $x,y\in K$.

Then
\[
\tau(x,y)=\left(\cj x,\phi'(x)y +\phi''(x)\cj y\right)
\quad\forall x,y\in K.
\]

The condition $\tau\circ\tau=id_{K^2}$ translates into:
\begin{gather}
\phi'(\cj x)\phi'(x)+\phi''(\cj x)\cj\phi''(x) = 1\ \ \forall x\in K \\
\phi'(\cj x)\phi''(x)+\phi''(\cj x)\cj\phi'(x) = 0 \ \ \forall x\in K
\end{gather}

Condition $(v)$ yields:
It follows that
\begin{gather}
\phi'(x)=\lambda^6\phi'(\lambda^2x)\ \ \forall x\in K, \lambda\in S\\
\phi''(x)= \phi''(\lambda^2x)\ \ \forall x\in K, \lambda\in S
\end{gather}

Condition $(iii)$ implies that $\phi'$ can be written as a finite
sum
\[
\phi'(x) = \sum_{k,l} a_{k,l} x^k\left(\cj x\right)^l \ \ (a_{k,l}\in K)
\]
Then
\[
\lambda^6\phi'(\lambda^2x) = 
\sum_{k,l} a_{k,l}\lambda^{6+2k-2l} x^k\left(\cj x\right)^l \ \ 
\forall x\in, \lambda\in S
\]
because $\cj\lambda=\lambda^{-1}$ for $\lambda\in S$.
Hence equation $(3)$ implies that $a_{k,l}=0$ unless $3+k=l$.
It follows that
\[
\phi'(x)=\left(\cj x\right)^3R(x\cj x)\ \forall x\in K
\]
for some polynomial $R\in K[T]$.

Similarily the equation $(4)$ implies that
\[
\phi''(x)=Q(x\cj x)\ \forall x\in K
\]
for some polynomial $Q\in K[T]$.

It follows that
\[
\phi''(\cj x)=\phi''(x)=Q(x\cj x)
\]
and
\[
\phi'(\cj x)=x^3R(x\cj x)
\]
for all $x\in K$.

Using these identities, equation $(1)$ transforms into
\begin{equation}
\left(x\cj x\right)^3\left(R(x\cj x)\right)^2
+Q(x\cj x)\cj Q(x\cj x)=1
\end{equation}

If 
\[
Q(z)=\sum_{k=0}^d c_kz^k
\]
with $c_d\ne 0$, then
the map
\[
z\mapsto Q(z)\cj Q(z)\quad (z\in k)
\]
is given by a polynomial $P\in k[T]$ of degree
$2d$ with top coefficient $c_d\cj c_d\ne 0$.

Hence equation $(5)$ implies
\[
z^3\left(R(z)\right)^2 + P(z)=1 
\]
for all $z\in k$ which can be written in the form $z=x\cj x=N_{K/k}(x)$.
The above equation thus holds for infinitely many $z\in k$ and therefore
implies an identity of polynomials.

However, the first term ($z^3\left(R(z)\right)^2$) is a polynomial
in $z$ of odd degree (unless $R\equiv 0$) 
while the second term ($P(z)$) is a polynomial
of even degree.

From this it follows that this identity can hold only in the degenerate
case where $R\equiv 0$ and $P\equiv 1$. 

Then $Q$ is a constant polynomial
taking an element of $S$ as value.
\end{proof}

\section{Proof of the theorem}
\begin{proof}
Due to the classification of $k$-forms of $G_m$ it follows that there
is a quadratic extension of fields $K/k$ with $\Gal(K/k)=\{id_K,\sigma\}$
such that $H^0$ arises
as described above in \S5.
We consider $G$ and $H$ and their actions $\mu$ resp.~$\mu_0$ on
$\A^4$
as defined in \S5 above.
We claim that the $H^0$-action $\mu_0$ on $E_0\simeq\A^4$ is not linearizable.
To see this, assume the contrary.
Recall that $E_0$ is a $H$-vector bundle over $\A^2$. Observe that
$(-1)\cj(-1)=1$ and that therefore $I=\{1,-1\}\subset S
=\{z\in K^*:z\cj z=1\}$. Regarded as subgroup of $H^0(K)\simeq S$,
the group $I$ acts on $E_0$ fixing precisely the zero-section.
Therefore we can employ the criterion of Masuda-Petrie (\S3)
 and deduce that
$E_0$ must be trivial as $H^0$-vector bundle over $\A^2$.
Hence, after applying a suitable $k$-automorphism of $E_0$, the
$H^0$-action on $E_0$ becomes linear.
Now 
$H$ is the semi-direct product of $\Z/2\Z=\{e,\tau\}$
and $H_0$ with
\begin{align*}
H_0(k)&=\{g\in G(K):\cj g=\left(g\right)^{-1}\}^0\\
&=S=\{\lambda\in G_m(K)=K^*: N_{K/k}(\lambda)=1\}.
\end{align*}
Linearity of the $H^0$-action on $E_0=k^4$ implies that there
is $H^0$-equivariant $k$-isomorphism between $E_0$ and the
tangent space at the fixed point $0\in E_0$.
Therefore we can identify $E_0=k^4$ with $K^2$ in such a way that
the $H^0$-action takes the form
\[
\lambda: (x,y)\mapsto (\lambda^2x,\lambda^3y)\ \ \forall \lambda\in S.
\]
Because $E_0$ is a $H$-vector bundle over $k^2$, conditions $(i)-(iii)$
of the proposition in \S6 are fulfilled. The group structure of
$H=\{e,\tau\}\ltimes H^0$ implies that conditions $(iv)$ and $(v)$
are fulfilled as well.

Hence we may invoke this proposition and deduce
that the action of $\tau$ on $E_0$ is already be given by
a $k$-linear map.
Since $H$ is generated by $\tau$ and $H^0$. it follows that
the $H$-action on $E_0$ is necessarily $k$-linear as soon as
the $H^0$-action is $k$-linear.
However, the action of $H$ is a $k$-form of the $G$-action on
$E$. It follows that the $G$-action on $E$ must be linearizable.
This is a contradiction, because the latter action is not linearizable
over $\C$.
\end{proof}

\end{document}